\documentclass[11pt,letterpaper]{amsart}
\usepackage{preamble}
\counterwithin{equation}{section}
\usepackage[section]{placeins}

\title{Multiplayer Games of War}

\author{Axel Adjei, Neil Krishnan, Elchanan Mossel}
\address{Department of Mathematics, Massachusetts Institute of Technology, Cambridge, MA 02139, USA}
\email{\{asadjei, neilk301, elmos\}@mit.edu}

\begin{document}



\begin{abstract}
A recent paper by Bhatia, Chin, Mani, and Mossel (2026) defined stochastic processes modeling the game of War for {\em two players} with $n$ cards. That paper showed that these models, assuming uniform random decks, are equivalent to the Gambler's Ruin problem and therefore have an expected termination time of $\Theta(n^2)$. In this paper, we generalize these models to {\em any number of players} $m$. We prove that the game with $m$ players is equivalent to a simple sticky random walk on an $(m-1)$-simplex; therefore, the termination time is the same as the absorption time of the simple sticky random walk.
Unlike Gambler's Ruin, this absorption time has not been analyzed before. We show that the absorption time of the walk and the termination time of the game are both $\Theta(n^2)$ for any number of players.
\end{abstract}

\maketitle
\section{Introduction}\label{sec:introduction}

\textit{War} is a simple card game. In each round of the game, all players place the top card of their hand face up. Whichever player plays the higher number wins the played cards, with a number of different variants used to determine ties. The game ends once one player collects all the cards, at which point this player is declared the winner. Recent work~\cite{bhatia2026will} introduced a number of mathematical models for this game with $2$ players and $n$ cards. In particular, \cite{bhatia2026will} showed that many variants of the $2$-player war game, where the initial decks are sampled uniformly and cards are played uniformly from the players' hands, are equivalent to the Gambler's ruin problem. 
This implies that the expected termination time of these variants is $\Theta(n^2)$.

We are interested in understanding mathematical models of the game of war when the number of players $m$ is arbitrary and for a large number of cards $n$. In our main result, we show that for cards drawn uniformly at random from the deck and under certain symmetry conditions, the process consisting of the decks of all players is equivalent to a \textit{sticky random walk} on the simplex. The term ``sticky'' refers to the fact that once the walk hits any lower-dimensional facet of the simplex, it stays there. 
Our results imply that the expected termination time of the game is the same as the absorption time of the sticky random walk at one of the vertices of the simplex. Interestingly, the absorption time of this walk seems not to have been analyzed before. We note that some problems involving random walks on a simplex are quite difficult; see e.g.~\cite{diaconis2021gambler,o20234}. 

For our problem of bounding the absorption time, we show it is $\Theta(n^2)$, where $n$ is the number of cards, for any number of players $m$, if $n/m$ is larger than a constant $\beta$ depending on the initial distribution of cards and other such parameters of the game.
This implies that the termination time of war game models is also $\Theta(n^2)$. 


The paper~\cite{bhatia2026will} also considers another variant called Martingale war or $f$-war where the chance of a card $a$ winning is proportional to $f(a)$ and where now cards are played from the top of the player's deck and returned to the bottom in arbitrary order. 
Under some conditions on $f$, it is shown in~\cite{bhatia2026will} that the termination time of the game is $\Theta(n^2)$. We prove similar results for any number of players. 


In the remainder of \cref{sec:introduction}, we will define the models for $m$-player war which we will consider in the rest of the paper. In \cref{sec:stickyrw}, we define the sticky random walk and show that these variants are equivalent to a sticky random walk on a simplex. In \cref{sec:absorptiontime}, using martingale arguments, we show that the sticky random walk terminates in $\Theta(n^2)$ time steps in expectation. This implies that the $m$ player variants also last $\Theta(n^2)$ turns in expectation. In \cref{sec:sim}, we use simulations to estimate the constant for the $\Theta(n^2).$ In \cref{sec:approxformula}, we motivate why finding an approximate formula for the expected time until some player loses all of their cards is interesting and propose an approximate formula for it, and in \cref{sec:furtherdirec}, we propose some further directions based on this formula.

We will now define terminology and models which will be referred to throughout the paper. Let the $m$ players be labeled $1, 2,\ldots, m$. We start with a deck $D$ with $n$ cards. Let $A_i(t)$ denote the cards in player $i$'s hand at time $t$. In each round, each player $i$ plays a card $a_i$ chosen through some process depending on the variant from their deck. If $A_i(t)$ is empty, then that player does not play a card, denoted by $a_i = \varnothing.$ Let $S_i = A_i(t) \setminus \{a_i\}$, i.e., $S_i$ is player $i$'s hand after playing $a_i$. Let a game tuple be the tuple $(a_1, \ldots, a_m, S_1, \ldots, S_m).$ We adapt the definition of a winning rule used in \cite{bhatia2026will}.

\begin{definition}
    A \textit{winning rule} $\mfp$ is a function which maps a game tuple $T = (a_1, \ldots, a_m, S_1, \ldots, S_m)$ to an $m$-dimensional vector $\mfp(T)$ where the $i$th entry, denoted $\mfp_i(T),$ is the probability of person $i$ winning. This map must satisfy the following conditions:
    \begin{itemize}
    \item If $a_i = \varnothing,$ then $\mfp_i(T) = 0.$
    \item For all permutations $\sigma$ of $[m],$
    $$\mfp_i(a_1,\ldots,a_m, S_1, \ldots, S_m) = \mfp_{\sigma(i)}(a_{\sigma(1)},\ldots,a_{\sigma(m)},S_{\sigma(1)},\ldots,S_{\sigma(m)}).$$
    \end{itemize}
\end{definition}

Additionally, a rule is said to be \textit{symmetric} if for all $i$ and all permutations $\sigma: [m] \to [m]$,
\begin{equation*}
    \mfp_i(a_1, \ldots, a_m, S_1, \ldots, S_m) = \mfp_{\sigma^{-1}(i)}(a_{\sigma(1)}, \ldots, a_{\sigma(m)}, S_1, \ldots, S_m).
\end{equation*}

Let us now define $\mfp$-war which is an $m$-player extension of \cite[Model 1]{bhatia2026will}.
\begin{model}[$\mfp$-war]\label{def:pwar}
$\mfp$-war with $m$ players can be modeled as a Markov chain on the state space $\Omega = \{(A_i)_{i\in [m]} : \bigsqcup A_i = D\}$ where the $A_i$ are disjoint subsets of $D$. A state represents the cards in each player's hand. We define $A_i(0)$ by choosing a vector $\x$ consisting of $m$ nonnegative integers whose sum is $n$ and then randomly shuffling the deck and partitioning the cards between the players such that $|A_i(0)| = x_i$. At the start of each turn, each player uniformly draws a card $a_i$ from their deck. From each state, with probability $\mfp_j(T)$ where $T$ is the game tuple, player $j$ wins the cards played while the other players lose the card they played. More formally,
$$A_i(t+1) = \begin{cases}
    A_j(t) \cup \{a_1, \ldots, a_m\} & i = j \\
    A_i(t) \setminus\{a_i\} & \text{otherwise}
\end{cases}.$$
\end{model}

\begin{remark}
    Note that the players who do not have cards do not contribute to the game. We keep them as they help separate similar states in the Markov chain. Thus, in this paper, we assume $A_i(t)$ is nonempty unless explicitly stated.
\end{remark}

Though \cite{bhatia2026will} does not mention the following variant, it is remarkably similar to actual war and can be dealt with using the same methods as for $\mfp$-war.

\begin{model}[Top card $\mfp$-war]\label{def:topcardpwar}
Top card $\mfp$-war is exactly the same as $\mfp$-war except players now draw from the top of their deck and when they win cards, the cards are randomly shuffled and added to the bottom of their deck.
\end{model}

In reality though, the cards might not be permuted before entering into the player's deck. For example, cards could be played in a pile with bottommost card of the pile being player 1's card and the topmost card being player $m$'s card. Then, when player $j$ adds the cards back to their deck, it is somewhat ordered, as there is a strong chance that the $j$th card from the bottom is ``better'' than the other cards as the $j$th player won the round. 

Therefore, we propose the following model which is the $m$-player extension of \cite[Model 2]{bhatia2026will} which allows for arbitrary ordering of the cards before adding them to the hand at the cost of more structure to the probability that a player wins a round. Each card is assigned a real number representing its strength. The probability that a given player wins a round is proportional to the strength of their card. The formal description of this model is as follows.

\begin{model}[$f$-war]\label{def:fwar}
Fix a function $f: \{1,\ldots, n\} \cup \{\varnothing\} \to \mathbb{R}^+$ such that $f(\varnothing) = 0$. We define $m$-player $f$-war to be a Markov chain on the space $\Omega = \{(A_i)_{i\in [m]} : \bigsqcup_{i = 1}^m A_i = [n]\}$. The $A_i$ are disjoint ordered tuples of elements in $[n]$.
We define $A_i(0)$ by choosing a vector $\x$ consisting of $m$ nonnegative integers whose sum is $n$ and then randomly shuffling the deck and partitioning the cards between the players such that $|A_i(0)| = x_i$.
The absorbing states of the Markov chain are those where there exists an $i$ such that $|A_i| = n$. At the start of each turn, player $i$ plays the card $a_i$ at the start of $A_i(t).$ With probability
$$\frac{f(a_i)}{\sum_{j = 1}^m f(a_j)},$$
player $i$ wins the played cards while the others lose the card they played. These new cards are added to the end of $A_i(t)\setminus\{a_i\}$ in arbitrary order.
\end{model}

We will also use the notation $A_{i, t}$ to denote $|A_i(t)|$ and $C_t$ to denote the set of players who have not run out of cards at time $t.$

The main result of the paper is that $\mfp$-war, top card $\mfp$-war, and $f$-war for certain $f$ all take $\Theta(n^2)$ turns, and we show this by first showing they are equivalent to a sticky random walk and then showing that the absorption time for sticky random walks is $\Theta(n^2).$


\section{The Sticky Random Walk}\label{sec:stickyrw}
In this section we first define the simple sticky random walk and show how it applies to $\mfp$-war and top card $\mfp$-war. We then define the fully general sticky random walk and show how it applies to $f$-war.

\subsection{Simple Sticky Random Walk}
We start by defining the simple sticky random walk.
\begin{definition}
    A \textit{simple sticky random walk} is a random walk in the state space $\Delta_n^{m-1} = \{(x_1, \ldots, x_m) \in \N^m: \sum_{i = 1}^m x_i = n\}$ which forms an $(m-1)$-simplex. Let $A_{i,t}$ be the $i$th coordinate of our state at time $t.$ Let $C_t$ be the set of $i$ where $A_{i,t}$ is positive at time $t$. From the state at time $t,$ we uniformly choose an index $j$ in $C_t.$ We then let $A_{j,t+1} = A_{j,t} + |C_t|-1$ and $A_{i,t+1} = A_{i,t} - 1$ for $i\neq j.$
\end{definition}

Let us now show that $\mfp$-war is equivalent to a simple sticky random walk which is analogous to \cite[Theorem 2.6]{bhatia2026will}.

\begin{proposition}\label{prop:equivtorandomwalk}
Consider a game of $\mfp$-war with a deck of $n$ cards where $\mfp$ is a symmetric winning rule. Assume that $\x$ is the distribution of the cards at the start of the game, i.e., $x_i$ is the number of cards in the $i$th players hand. Then, $(A_{i, t})_{i\in [m]}$ is a simple sticky random walk starting at $\x$.
\end{proposition}
\begin{proof}
Recall that in the definition of $\mfp$-war, $A_i(0)$ is defined by shuffling the deck randomly and partitioning the cards so that the $i$th player gets $x_i$ cards. Thus, given that the size of $A_i(0)$ is $x_i$ for some vector $\x \in \Delta_n^{m-1},$ we see $(A_i(0))_{i\in [m]}$ is uniformly distributed among all sequences of hands where each player $i$ has hand size $x_i$ cards.

Let $t$ be a time at which the game has not ended yet. Assume that $(A_i(t))_{i\in [m]}$ are uniformly distributed in partitions of the deck where the size of the $i$th partition is $A_{i,t}$.
It follows that the tuple $(a_i)_{i \in C_t}$ is uniform across all tuples of cards of size $|C_t|$. Let $|C_t| = r$ and, without loss of generality, assume that $C_t = [r].$ If we just fix the set of cards played to be $B = \{b_i: i \in [r]\}$ and each $S_i,$ the tuple $(a_i)_{i \in C_t}$ is equally likely to be any ordering of $B$. Because the winning rule is symmetric and the remaining decks of each of the players is fixed, the probability of player $i$ winning by playing a card $b$ is equal to the probability of player $j$ winning if they play $b$, so there is some fixed $p_{b_i}$ which is the probability that the person placing $b_i$ wins. Note that $\sum_{ i = 1}^r p_{b_i} = 1.$ Thus, the probability player $i$ wins is
$$\frac{1}{r!}\sum_{\sigma \in \fS_r} \mfp_i(b_{\sigma(1)},\ldots,b_{\sigma(m)},S_1,\cdots,S_m) = \frac{1}{r!}\sum_{\sigma \in \fS_r} p_{b_{\sigma(i)}} = \frac{1}{r} \sum_{ i = 1}^r p_{b_i} = \frac{1}{r}.$$
This implies that $(A_{i, t})_{i\in [m]}$ is a sticky random walk assuming that $(A_i(t))_{i \in [m]}$ are always uniformly distributed in partitions where the size of the $i$th partition is $A_{i,t}$.

We just need to show that the $A_i(t)$ remain uniformly distributed now. Using induction, it suffices to show that if $(A_{i}(t))_{i \in [m]}$ are uniformly distributed in partitions where player $i$ has $A_{i,t}$ cards, then given that player $j$ wins, $(A_i(t+1))_{i\in [m]}$ is uniformly distributed in sequences where player $i$ has $A_{i,t+1}$ cards. We do this by describing the procedure of a turn in a different way. We first reveal the cards played $B$ uniformly chosen in sets of size $|C_t|.$ Note that given these cards, as we saw earlier, each player wins with equal probability. Assume that player $j$ wins. Reveal the set $S_j$ uniformly distributed in subsets of $D\setminus B$ of size $A_{j,t}-1$ Finally, reveal the remaining sets $(S_i)_{i \in [m]\setminus\{j\}}$ which is uniformly distributed in partitions of $D\setminus (B \cup S_j)$ where the part indexed by $i$ has size $A_{i,t}-1.$ Player $j$ wins, so their hand is now $A_{j}(t+1) = S_j \cup B,$ which is uniformly distributed in hands of size $\lvert S_j\cup B \rvert.$ Because $(S_i)_{i \in [m]\setminus\{j\}}$ are uniformly distributed in partitions of appropriate size of the remaining cards, $A_i(t+1) = S_i$ for $i\neq j$ maintains that property.
\end{proof}

\begin{remark}
Note that by the same proof as \cref{prop:equivtorandomwalk}, top card $\mfp$-war can also be modeled by the simple sticky random walk. The critical step is that the cards drawn from the top of each player's deck are uniform in all sets of $|C_t|$ cards and the deck of the player who wins is uniform in ordered tuples not containing the played cards, so the hand formed by combining these two stacks of cards must also be uniform in all ordered tuples of that size. Because the proof is so similar, we omit it.
\end{remark}

In \cref{sec:absorptiontime}, we prove \cref{thm:srwabsorptiontime} which we state in the next subsection. A corollary of \cref{thm:srwabsorptiontime} relevant to simple sticky random walks is the following.
\begin{corollary}\label{cor:simplesrw}
    Let $n$ and $m\geq 2$ be positive integers. Assume that there is some $\alpha < 1$ such that for all $i\in [m],$ $A_{i,0}\leq \alpha n.$ Assume that $n/m$ is larger than some constant $\beta$ depending on $\alpha$. Then the expected absorption time of the simple sticky random walk is $\Theta(n^2).$
\end{corollary}
This corollary implies that the expected length of $\mfp$-war and top card $\mfp$-war is $\Theta(n^2)$
\begin{corollary}\label{cor:pwarlength}
    Let $n$ and $m\geq 2$ be positive integers. Let $\alpha < 1$ be a fixed constant. Consider a game of $\mfp$-war or top card $\mfp$-war with $n$ cards and $m$ players such that in the starting state $\x,$ we have that $x_i \leq \alpha n$ for all $i$. Assume that $n/m$ is larger than some constant $\beta$ depending on $\alpha$. The expected length of the game is then $\Theta(n^2).$
\end{corollary}

\subsection{Sticky Random Walk}

We now define the sticky random walk. Note that when we say a martingale $M_t$ increases or decreases, we mean that $M_{t+1} > M_t$ or $M_{t+1} < M_t,$ respectively.

\begin{definition}\label{def:srw}
    A \textit{sticky random walk} is a random walk $(A_{i,t})_{i \in [m]}$ on $\{(x_1, \ldots, x_m)\in \R^m: \forall x_i, x_i \geq 0, \sum_{ i = 1}^m x_i = n\}$ that satisfies the following properties for all $i$. In what follows, we use $r$ to denote the number of $A_{i,t}$ which are nonzero at time $t.$
    \begin{itemize}
    \item $(A_{i,t})_{t \in \N}$ is a martingale.
    \item If $A_{i,t} = 0,$ then $A_{i,t+1} = 0.$ 
    \item If $A_{i,t}$ increases, $A_{j,t}$ decreases for all $j\neq i.$ 
    \item The probability that $A_{i,t}$ increases is between $p_- = p_-(r) = \Omega(1/r)$ and $p_+ = p_+(r) = O(1/r)$ for some functions $p_-(r)$ and $p_+(r).$
    \item For some real numbers $a_-,a_+ > 0,$ if $A_{i,t}$
    decreases, we have $a_- \leq A_{i,t}-A_{i,t+1} \leq a_+.$
    \item For some real number $n,$ $\sum_{i = 1}^m A_{i,t} = n$
\end{itemize}
    States in the sticky random walk can be thought of as tuples $(A_{i,t})_{i \in [m]}$ and the absorbing state is when all of the $A_{i,t}$ are $0$ except for one.
\end{definition}

Notice that the simple sticky random walk is a sticky random walk where $p_- = p_+ = a_- = a_+ = 1.$ In \cref{sec:absorptiontime}, we will prove the following result, showing the sticky random walk absorption time is $\Theta(n^2)$.
\begin{theorem}\label{thm:srwabsorptiontime}
    Let $n$ be a positive real number and $m\geq 2.$ Assume that there is some $\alpha < 1$ such that for all $i \in [m]$ $A_{i,0} \leq \alpha n.$  Assume that $n/m$ is larger than some constant $\beta$ depending on $\alpha, a_-,a_+,p_-,p_+$. Then, the expected absorption time of a sticky random walk is $\Theta(n^2).$
\end{theorem}

This general sticky random walk can model $f$-war for certain $f.$ Note that in $f$-war, we also use $n$ to denote the total number of cards, so for the remainder of this section, we will use $n'$ to denote the total number of cards instead.
\begin{proposition}\label{prop:f}
    Let $f:\{1,\ldots,n'\}\cup\{\varnothing\}\to\R^+$ be a function where $f(\varnothing) = 0,$ and there exists $p_-(r) = \Omega(1/r)$ and $p_+(r) = O(1/r)$ so that for every set $S \subseteq [n']$ where $|S| = r$ and element $a \in S,$
    $$p_-(r) \leq \frac{f(a)}{\sum_{k \in S} f(k)} \leq p_+(r),$$
    and $ \max_{i \in [n']} f(i) \leq cx$ for some constant $c$ where $x = \min_{i \in [n']} f(i).$ Then, the random walk defined by
    $$\left(\frac{1}{x} \sum_{k \in A_i(t)} f(k)\right)_{i \in [m]},$$
    is a sticky random walk.
\end{proposition}
\begin{proof}
    Let $A_{i,t} = \frac{1}{x} \sum_{j \in A_i(t)} f(j).$ We organize this proof through the following list which proves the corresponding properties listed in \cref{def:srw}.
    \begin{itemize}
    \item We first must show that each component is a martingale. Assume without loss of generality that only players $1$ through $r$ still have cards. Let $b_j$ for $1\leq j \leq r$ be the top card of player $j$'s deck. We then have that
    $$\EE\left[A_{i,t+1}\right] = \sum_{k \in A_i(t+1)} f(k) + \frac{b_i}{\sum_{j = 1}^r b_i} \left(\sum_{\substack{j = 1 \\ j\neq i}}^r b_i \right) + \frac{\sum_{\substack{j = 1 \\ j\neq i}}^r b_i}{\sum_{j = 1}^r b_i} (-b_i) = \sum_{k \in A_i(t+1)} f(k) = A_{i,t},$$
    so each coordinate is a martingale. 
    \item If $A_{i,t} = 0,$ then player $i$ has no cards, so $A_{i,t+1}$ stays $0.$
    \item If $A_{i,t}$ increases, then player $i$ gained cards which means that all other players lost cards, so all other $A_{j,t}$ for $j\neq i$ decreases.
    \item The probability of $A_{i,t}$ increasing is $$\frac{f(a_i)}{\sum_{j = 1}^m f(a_j)} = \frac{f(a_i)}{\sum_{j = 1}^r f(a_j)},$$
    which is between $p_-(r)$ and $p_+(r)$ as given by the proposition statement.
    \item If $xA_{i,t}$ decreases, it decreases by at least $\min_{i \in [n']} f(i)$ and at most $\max_{i \in [n']} f(i),$ so the $i$th coordinate which is $A_{i,t}$ decreases by at least $a_- = 1$  and at most $a_+ = \frac{1}{x} \max_{i \in [n']} f(i).$ As given by the proposition statement, $a_-= 1$ and $a_+ \leq c$.
    \item We have that
    $$\sum_{i = 1}^m A_{i,t} = \frac{1}{x} \sum_{k = 1}^{n'} f(k),$$
    which is a constant.
    \end{itemize}
\end{proof}
\begin{example}
    Let $f(i) = n'+i$ for a card $i$ which is the $f$ mentioned in \cite[Claim 3.8]{bhatia2026will}. As always $f(\varnothing) = 0.$ This $f$ satisfies all of the properties in the statement of \cref{prop:f}. For a given $r,$ the probability of player $i$ winning is minimized when $a_i = 1$ and $a_j = n'$ for $j\neq i$, so the probability is at least $\frac{n'+1}{n'+1+(r-1)(2n')} = \Omega(1/r) = p_-(r),$ and the probability of player $i$ winning is maximized when $a_i = n'$ and $a_j = 1$ for $j\neq i,$ so the probability is at most $\frac{2n'}{2n'+(r-1)(n'+1)} = O(1/r) = p_+(r).$ Finally $\min_{i \in [n']} f(i) = n'+1$ which is a constant factor off of $\max_{i \in [n']} f(i) = 2n'.$ Thus \cref{prop:f} applies to $f$-war for this choice of $f.$
\end{example}

Using \cref{thm:srwabsorptiontime}, we can show that the expected length of $f$-war for the $f$ specified in \cref{prop:f} is $\Theta({n'}^2).$
\begin{corollary}
    Let $n'$ and $m\geq 2$ be positive integers. Let $\alpha < 1$ be a fixed constant. Consider a game of $f$-war with $n'$ cards and $m$ players such that in the starting state $\x,$ we have that $x_i \leq \alpha n$ for all $i.$ Also assume that $f$ satisfies the conditions in the statement of \cref{prop:f}. Then the expected length of $f$-war is $\Theta({n'}^2).$
\end{corollary}

\section{Bounds on Absorption Time}\label{sec:absorptiontime}
We will now prove \cref{thm:srwabsorptiontime}. We split the proof into two parts. We first establish propositions upper bounding the expected amount of time where the number of positive martingales is $r,$ and in the next subsection, we use a super and sub martingale to show the absorption time is $\Theta(n^2).$

\subsection{Expected Amount of Time with $r$ Positive Martingales}
In this subsection, we determine upper bounds on the expected amount of time where the number of positive martingales is $r\leq m$.
We first consider when $n\leq r^2.$

\begin{proposition}\label{prop:smallntime}
    If $r \leq n\leq r^2,$ the expected amount of time when the number of positive martingales is $r$ is $O(n/r).$
\end{proposition}
\begin{proof}
We will examine the random walk of $M_t = \min \{A_{i,t}: i \in C_t\}$. At the start, note $M_0 \leq q$ where $q = n/r.$ If $M_{t-1} = b,$ we know $M_t \leq b-a_-$ with probability at least $1-p_+$. Otherwise, with probability at most $p_+,$ we have that $M_t$ will increase by some amount. Note that $M_t \leq q$ always. When $M_t = 0,$ the number of positive martingales will be less than $r.$ Therefore, we want to find the expected amount of time until $M_t$ hits $0.$

Consider the random walk $Z_t$ where $Z_0 = q,$ $Z_t = Z_{t-1} - a_-$ with probability $1-p_+,$ and $Z_t = q$ with probability $p_+.$ It is clear that the expected time when $Z_t$ hits $0$ is larger than the expected time when $M_t$ hits $0.$

Let us now find the expected time when $Z_t$ hits $0.$ Assume at time $s,$ $Z_s = q.$ The probability that $Z_t$ for $t\geq s$ continually decreases until hitting $0$ is least
$$\left( 1-p_+\right)^{q/a_-} \geq \left( 1 - p_+\right)^{r/a_-}\geq p^*,$$
for some fixed $p^*$ which varies only with $p_+$ and $a_-.$ Critically, we use the fact that $q \leq r$ as $n\leq r^2$ here. If $Z_t$ does not decrease until hitting $0$, for $s+1 \leq t \leq s+q/a_-,$ we have that $Z_t = q.$ The largest possible gap then between consecutive times when $Z_t = q$ is $q/a_-.$ Thus, the expected time until $Z_t = 0$ is at most
$$p^* \left( \frac{q}{a_-} \right) + p^*(1-p^*) \left( 2\frac{q}{a_-} \right) + p^*(1-p^*)^2 \left(3 \frac{q}{a_-} \right) + \cdots = O(q) = O\left(\frac{n}{r}  \right).$$
\end{proof}

Now consider if $n\geq r^2.$ Let $\x$ be an $r$-dimensional vector with the $i$th dimension being the $i$th positive martingale. Without loss of generality, these martingales are $A_{1,t}$ through person $A_{r,t}.$ Let $E(\x)$ denote the expected amount of time until one of the martingales hits $0$. Let $\v_i$ be the vector representing how $\x$ will change assuming the $i$th player wins. Let $b(\x)$ denote the minimum of the entries in $\x.$ Let $p_i$ be the probability the $i$th martingale increases.

\begin{proposition}\label{prop:largentime}
    If $n\geq r^2,$ the expected amount of time when the number of positive martingales is $r$ is $O(n^2/r^3).$
\end{proposition}
\begin{proof}
    We have that $E(\x) = 0$ if $b(\x) = 0.$ Otherwise, we have
    $$E(\x) = 1 + \sum_{i = 1}^r p_i E(\x+\v_i).$$
    Consider a function $h(\x)$ such that $h(\x) = 0$ if $b(\x) = 0$ and otherwise
    $$h(\x) \geq 1 + \sum_{i = 1}^r p_i h(\x+\v_i).$$
    Let us show that $E(\x) \leq h(\x).$ Let $u(\x) = h(\x) - E(\x)$ so that $u(\x) = 0$ if $b(\x) = 0$ and otherwise
    $$u(\x) \geq \sum_{i = 1}^r p_i u(\x+\v_i).$$
    Let $\y$ minimize $u(\y).$ Assume $u(\y) < 0.$ Note then $b(\y) \neq 0.$ Then 
    $$\sum_{i = 1}^m p_i u(\y+\v_i) \leq u(\y),$$
    so $u(\y+\v_i) = u(\y)$ for all $i,$ and this applies for all $\y$ achieving this minimum. Consider a particular $i \in [m].$ Note assuming that starting at $\y,$ we keep letting the $i$th martingale increase. The value of $u$ will stay the same. Eventually, one of the martingales will hit $0$ at which time the value of $u$ will be $0.$ This is a contradiction.

    We now need to find such a function $h(\x).$ Consider
    $$h(\x) = \frac{1}{p_- a_-^2(r-1)^2} \left( \left( \frac{n}{r} \right)^2 - \left( \frac{n}{r} - b(\x) \right)^2 \right) + \frac{1}{p_- a_-^2(r-1)^2}\cdot 2a_+ \left(\frac{1}{p_-} - 1\right) \cdot b(\x).$$
    Note that if $b(\x) = 0,$ then $h(\x) = 0.$ Without loss of generality, assume that $A_{1,t}$ achieves the minimum value of $b(\x).$ Let $a_i \geq a_-$ for $2\leq i \leq r$ be the value by which $b(\x)$ decreases if $A_{i,t}$ increases and let $\ell$ be the value by which $b(\x)$ changes if $A_{1,t}$ increases. Note that $\ell$ can be negative. Notice that in $h(\x)$ the only components that change with respect to $\x$ are
    $$2\frac{1}{p_- a_-^2(r-1)^2}\left( \frac{n}{r} + a_+ \frac{1 - p_-}{p_-} \right) b(\x) - \frac{1}{p_- a_-^2(r-1)^2}b(\x)^2.$$
    Note that
    $$\sum_{ i = 1}^r p_i b(\x+\v_i) = b(\x) +  p_1 \ell - \sum_{i = 2}^r p_i a_i,$$
    and
    $$\sum_{i = 1}^r p_i b(\x+\v_i)^2 \geq p_1(b(\x)+\ell)^2 + \sum_{i = 2}^r p_i(b(\x)-a_i)^2 = b(\x)^2 + 2\left(p_1\ell - \sum_{i = 2}^r p_i a_i\right)b(\x) + p_1\ell^2 + \sum_{ i = 2}^r p_i a_i^2.$$
    Thus,
    \begin{align*}
        h(\x) - \sum_{i = 1}^r p_i h(\x+\v_i) = &\frac{2}{p_- a_-^2(r-1)^2}\left( \frac{n}{r} + a_+ \frac{1 - p_-}{p_-} \right) \left(b(\x) - \left(b(\x) + p_1\ell - \sum_{i = 2}^r p_i a_i  \right)  \right),\\
        - &\frac{1}{p_- a_-^2(r-1)^2} \left(b(\x)^2 - \left( b(\x)^2 + 2\left(p_1\ell - \sum_{i = 2}^r p_i a_i\right)b(\x) + p_1\ell^2 + \sum_{ i = 2}^r p_i a_i^2 \right)  \right),\\
        = &\frac{1}{p_- a_-^2(r-1)^2}\left(2\left( \frac{n}{r} + a_+ \frac{1 - p_-}{p_-}-b(\x) \right)\left( -p_1 \ell + \sum_{i = 2}^r p_i a_i \right)+p_1\ell^2+\sum_{i = 2}^r p_ia_i^2\right).
    \end{align*}
    Let us show derivative with respect to $\ell$ is negative. It suffices to just compute the derivative of expression inside the first set of parenthesis. The derivative of the left term is
    $$-2p_1\left( \frac{n}{r} + a_+ \frac{1 - p_-}{p_-}-b(\x) \right) \leq -2p_1 a_+ \frac{1 - p_-}{p_-} \leq -2p_1 \frac{-p_1}{p_1}a_+ = -2(1-p_1)a_+,$$
    while the derivative of the right term is
    $$2\ell p_1 \leq 2\sum_{i = 2}^r p_i a_i \leq 2(1-p_1)a_+.$$
    Thus, because the derivative with respect to $\ell$ is negative, $h(\x) - \sum_{i = 1}^r p_i h(\x+\v_i)$ is minimized when $\ell$ is as large as possible, i.e., when $p_1\ell = \sum_{i = 2}^r p_i a_i.$ In that case, we have
    \begin{align*}
    h(\x) - \sum_{i = 1}^r p_i h(\x+\v_i) &= \frac{1}{p_- a_-^2(r-1)^2}\left(p_1\ell^2+\sum_{i = 2}^r p_ia_i^2  \right) \geq \frac{1}{p_- a_-^2(r-1)^2} \cdot p_1 \ell^2,\\
    &\geq \frac{1}{p_- a_-^2(r-1)^2} p_1((r-1)a_-)^2 \geq 1.
    \end{align*}
    Therefore,
    $$h(\x) \geq 1 + \frac{1}{r} \sum_{i = 1}^r p_i h(\x+\v_i),$$
    so $E(\x) \leq h(\x) = O(n^2/r^3).$
\end{proof}

\subsection{Expected absorption time}

Now we can determine the expected absorption time of the sticky random walk. We start by introducing a submartingale $(X_t)_{t\in\N}$ and a supermartingale $(Y_t)_{t\in\N}.$ In the following definition, $r$ is the number of positive martingales at time $s.$

\begin{lemma}\label{lem:cubesmartingale}
    Let $$X_t' = \sum_{i = 1}^m A_{i,t}^3 - \sum_{s = 0}^{t-1} w(n,r),$$ where $w(n,r) = \Theta(nr+ r^3).$ Then for a certain choice $w(n,r) = y(n,r)$, $(X_t)_{t \in \mathbb{N}} = (X_t')_{t\in \mathbb{N}}$ is a submartingale and for another such choice $w(n,r) = v(n,r)$ $(Y_t)_{t \in \mathbb{N}} = (X_t')_{t\in \mathbb{N}}$ is a supermartingale. Furthermore $y(n,r) = \Omega(n).$
\end{lemma}

\begin{proof}
    Consider the filtration $\FF_t = \sigma(A_{j,s}: j \in [m], s\leq t).$ We need to show that $\EE[X_{t} \mid \FF_{t-1}] \leq X_{t-1}$ and $\EE[Y_t\mid\FF_{t-1}] \geq Y_{t-1}.$ Let $\nu_{ij}$ for $i,j \in C_{t-1}$ be the amount $A_{i,t-1}$ changes to get to $A_{i,t}$ assuming the $A_{j,t}$ increases. Without loss of generality, assume that only martingales $A_{1,t-1}$ through $A_{r,t-1}$ are still positive at time $t-1.$ We then have
    \begin{align*}
    \EE\left[\sum_{i = 1}^r A_{i,t}^3 \bigr\vert \FF_{t-1}\right] &= \sum_{j = 1}^r p_j \sum_{i = 1}^r (A_{i,t-1} +\nu_{ij})^3,\\ 
    &= \sum_{j=1}^r p_j\left(  \sum_{i = 1}^r A_{i,t-1}^3 + 3 \sum_{i = 1}^r A_{i,t-1}^2 \nu_{ij} + 3\sum_{i = 1}^r A_{i,t-1} \nu_{ij}^2 + \sum_{i =1}^r \nu_{ij}^3\right),\\
    &= \sum_{i = 1}^r A_{i,t-1}^3 +3\sum_{i = 1}^r A_{i,t-1}^2 \left(\sum_{j = 1}^r p_j \nu_{ij}\right) + 3\sum_{i = 1}^r A_{i,t-1} \left( \sum_{j = 1}^r p_j\nu_{ij}^2\right) + \sum_{j = 1}^r p_j \sum_{i =1}^r \nu_{ij}^3,\\
    &= \sum_{i = 1}^r A_{i,t-1}^3 + 3\sum_{i = 1}^r A_{i,t-1} \left( \sum_{j = 1}^r p_j\nu_{ij}^2\right) + \sum_{j = 1}^r p_j \sum_{i =1}^r \nu_{ij}^3.
\end{align*}
For any $i$,
$$\Theta_-(r) = p_- \left((r-1)a_-\right)^2 + \left( 1 - p_+\right)\left( a_-^2\right) \leq \sum_{j = 1}^r p_j \nu_{ij}^2 \leq p_+ \left((r-1)a_+\right)^2 + \left( 1 - p_-\right)\left( {a_+}^2\right) = \Theta_+(r),$$
where we use $\Theta_-$ and $\Theta_+$ to differentiate the two asymptotic expressions. We also have that for any $j,$
$$\Theta_-(r^3) = (a_-(r-1))^3 - (r-1){a_+}^3\leq \sum_{i = 1}^r \nu_{ij}^3 \leq (a_+(r-1))^3 -(r-1)a_-^3 = \Theta_+(r^3).$$
Let
$$v(n,r) = 3\sum_{i = 1}^r A_{i,t-1} \Theta_+(r) + \sum_{j = 1}^r p_j \Theta_+(r^3) = n\Theta_+(r) + \Theta_+(r^3) = \Theta(nr+r^3),$$
and
$$y(n,r) = 3\sum_{i = 1}^r A_{i,t-1} \Theta_-(r) + \sum_{j = 1}^r p_j \Theta_-(r^3) = n\Theta_-(r) + \Theta_-(r^3) = \Theta(nr+r^3).$$
Thus,
$$\sum_{i = 1}^r A_{i,t-1}^3 + y(n,r) \leq \EE\left[\sum_{i = 1}^r A_{i,t}^3 \bigr\vert \FF_{t-1}\right] \leq \sum_{i = 1}^r A_{i,t-1}^3 + v(n,r).$$
Note that $\Theta_-(r)$ has a positive coefficient for its order $r$ term. Therefore, setting, $r = 2,$ we see that $y(n,2) = \Theta(n)$ and this is a lower bound for $y(n,r)$ modulo a constant.

Observe that $w(n,r)$ is precisely the maximum and minimum change in $\sum_{ i = 1}^r A_{i,s}^3$ from $s = t-1$ to $s = t$ in expectation in both the supermartingale and submartingale cases. Thus,
$$\sum_{s = 0}^{t-1} w(n,r),$$
is the cumulative change in $\sum_{i = 1}^m A_{i,t}^3$ which makes $X_t'$ a supermartingale if $w = v$ and a submartingale if $w = y.$

\end{proof}

We are now ready to prove that the expected length of a sticky random walk is $\Theta(n^2).$
\begin{proof}[Proof of \cref{thm:srwabsorptiontime}]
    Let $\tau$ be the time when only one martingale is positive. Note that $\EE[\tau] < \infty$ from standard Markov chain analysis. Substituting the event of the game ending with the event of $A_{1,t}$ increasing at every time step until the game ends, we have a geometric random variable as $A_{1,t}$ either increases until the game ends or it does not and we, at worst, reset with $A_{1,t}$ being very small again. This random variable must have a higher expectation than $\tau$, but the expectation is finite as the probability of $A_{1,t}$ increasing until the game ends is lower bounded. Moreover, the differences $|X_{t+1}'-X_t'|$ are bounded above by $n^3 + v(n,m)$ as $n^3$ bounds how much the first sum of $X_t'$ can change and $v(n,m)$ bounds how much the second sum of $X_t'$ can change. Thus, we can apply the optional stopping theorem where we stop the random walk at $t = \tau.$

    Let us consider the submartingale from \cref{lem:cubesmartingale} first. We will use this to show that $\EE[\tau] = O(n^2).$ By the optional stopping theorem, we have
    $$\EE\left[\sum_{ i = 1}^m A_{i,0}^3\right] = \EE[X_0] \leq \EE[X_\tau] = n^3 - \EE\left[ \sum_{s = 0}^{\tau-1} y(n,r) \right],$$
    so
    $$\EE\left[ \sum_{s = 0}^{\tau-1} y(n,r) \right] \leq n^3 - \EE\left[\sum_{ i = 1}^m A_{i,0}^3\right] \leq n^3.$$
    Thus, we have
    $$n^3\geq \EE\left[ \sum_{s = 0}^{\tau-1} y(n,r) \right] = \EE[\tau] \Omega(n) \Rightarrow \EE[\tau] = O(n^2).$$
    
    We now need to show that $\EE[\tau] = \Omega(n^2).$ Recall we have $A_{i,0}\leq \alpha n.$ We only lose a constant factor. If we take the supermartingale from \cref{lem:cubesmartingale} and apply the optional stopping theorem, we get
    $$
        \EE\left[\sum_{ i = 1}^m A_{i,0}^3\right] = \EE[Y_0] \geq \EE[Y_{\tau}] = n^3 - \EE\left[\sum_{s = 0}^{\tau - 1}v(n,r)\right].
    $$
    Letting $T_r$ be the random variable for the amount of time spent with $r$ players in the game, we have
    $$\EE\left[\sum_{ i = 1}^m A_{i,0}^3\right] \geq n^3 - \sum_{r = 2}^{m}v(n,r)\EE[T_r].$$
    Thus,
    $$\sum_{r = 2}^{m}v(n,r)\EE[T_r] \geq n^3 - \EE\left[\sum_{ i = 1}^m A_{i,0}^3\right] \geq cn^3,$$
    for some $c$ because $A_{i,0} \leq \alpha n$ for some $\alpha < 1.$

    Note that if we copy the strategy from the upper bound, we get that
    $$v(n,m) \EE[\tau] \geq \sum_{r=2}^m v(n,r) \EE[T_r] \geq cn^3,$$
    so $\EE[\tau] \geq cn^3/v(n,m).$ For constant $m,$ note that $v(n,m) = O(n),$ so $\EE[\tau] = \Omega(n^2).$ Thus, it suffices to consider only sufficiently large $m.$
    
    We will now proceed with cases based on how large $n$ is compared to $m.$

    \textbf{Case 1}: Assume $n\geq m^2.$ Recall from \cref{prop:largentime}, we have $\EE[T_r] = O(n^2/r^3).$ Assuming no other constraints on $\EE[T_r],$ let us minimize $\sum_{r = 2}^m \EE[T_r],$ the expected length of the game. The coefficient of $\EE[T_r]$ decreases as $r$ decreases, so for the same increase in $\EE[T_r],$ the increase contributes less to the inequality if $r$ is smaller. Thus, we want to maximize $\EE[T_r]$ where the coefficient of $\EE[T_r]$ is larger. The largest possible value of $\EE[T_r]$ is $O(n^2/r^3),$ and
       $$\sum_{r=2}^m v(n,r)\EE[T_r]
        \leq \sum_{r=2}^\infty \Theta(nr+r^3) \cdot O\left( \frac{n^2}{r^3} \right) \leq \sum_{r=2}^\infty \Theta(nr) \cdot O\left( \frac{n^2}{r^3} \right)
        \leq \kappa n^3,$$
    for some $\kappa$ which generally is larger than $c.$ This implies that truncating the sum from $r = r'$ for fixed $r' > 2$ will make it less than $cn^3.$ Thus, we need $\EE[T_r]$ to be at its maximum value for at least $r \geq r'.$ Thus, the expected length of the game is at least some constant times $n^2/{r'}^3$ which is $\Omega(n^2).$

    \textbf{Case 2}: Now consider $m^{5/3} \leq n \leq m^2.$ We use identical reasoning to the $n\geq m^2$ case except we must now use \cref{prop:smallntime} as well to bound $\EE[T_r].$ We have
    \begin{align*}
        &\sum_{r = 2}^m v(n,r)\EE[T_r],\\
        \leq & \sum_{r = 2}^{\lfloor \sqrt{n}\rfloor} v(n,r) O\left( \frac{n^2}{r^3} \right) + \sum_{r = \lfloor \sqrt{n}\rfloor+1}^{m} v(n,r) O\left( \frac{n}{r} \right),\\
        \leq &\sum_{r = 2}^{\infty} \Theta(nr) O\left( \frac{n^2}{r^3} \right) + \sum_{r = 1}^{m} (O(n^2) + O(nr^2)),\\
        \leq &\kappa n^3 + O(mn^2) + O(nm^3),\\
        \leq &\kappa n^3 + O(n^{13/5}) + O(n^{14/5}),
    \end{align*}
    where $\kappa > c.$ For large $m$ and $n,$ we have that $O(n^{13/5})$ and $O(n^{14/5})$ are insignificant compared to $\kappa n^3.$ Again, we can truncate the sum from $r = r'$ to $m$ in order to shrink $\kappa$ where $r'$ is a constant with respect to $n$ and $m,$ so the expected length of the game will be $\Omega(n^2)$ again.

    \textbf{Case 3}: Assume that $\beta m \leq n \leq m^{5/3}$ where $\beta$ is a sufficiently large constant not depending on $n$ or $m.$ Let $q$ be the number of time steps needed until the number of positive martingales remaining is $n^{3/5}.$ Consider the first $\sqrt{n}$ turns. Any martingale which is still positive after these turns must have size at least $a_-\sqrt{n}$ or have increased during at least one of these turns. There are at most $\frac{\sqrt{n}}{a_-}$ martingales with size at least $a_-\sqrt{n}$ and at most $\sqrt{n}$ martingales which increase during the first $\sqrt{n}$ turns, so the number of martingales remaining is $O(\sqrt{n}) < n^{3/5}.$ Thus $q\leq \sqrt{n}.$ Let us first show that with nonzero probability, no martingale increases $\gamma n/m$ or more times in the first $q$ turns where $\gamma$ is a sufficiently small constant such that $\alpha n +\gamma a_+(m-1) < \alpha'n$ for some $\alpha' < 1.$ This probability is at most, by a Markov bound, the expected number of people who win $\gamma n/m$ or more times which is at most
    $$m \dbinom{q}{\gamma n/m}p_+(n^{3/5})^{\gamma n/m} = O\left( m \frac{(\sqrt{n})^{\gamma n/m}}{(n^{3/5})^{\gamma n/m}} \right),$$
    and for sufficiently large $n$ and $m$ as well as sufficiently large $n/m = \beta \approx 10/\gamma,$ this fraction is at most $1/2.$ Therefore, after the first $q$ turns, with probability at least $1/2,$ no one wins $\gamma n/m$ or more times. Therefore, there are $n^{3/5}$ positive martingales left each of which is at most $\alpha n + \gamma a_+(m-1) < \alpha' n$ where $\alpha' < 1$. Imagine that we are starting the game now with this distribution of cards. Our new $m$ is $n^{3/5},$ so we are now in the case of $n\geq m^{5/3}.$ That argument finishes this case.

    
\end{proof}

\begin{remark}\label{rem:n/mlargerthanconstant}
Note that the constraint of $n/m$ being larger than some constant is necessary. Consider the simple sticky random walk where the cards are divided so that the absolute difference between the number of cards any two players have is at most $1$. If $n = 2m,$ then after the first two turns, with high likelihood, there are two players each of whom has roughly $m$ cards. The sticky random walk becomes a Gambler's ruin problem, and it is known that the expected length of the game is the product of the number of cards each of the players has which is $m^2 = \Theta(n^2).$ In fact, this logic holds for $n = cm$ where $c > 1$ and shows that the number of turns is $\Theta(n^2).$ On the other hand, if $n = m$ then there is exactly one turn. Thus, a sufficiently large constant $c$ is necessary.

There is a transition from $n = m$ to $n = cm$ for $c > 1$ where the number of turns is $\Theta(n^2).$ If $m \leq n = m(1+o(1)),$ notice that there are $n-m = o(m)$ players with $2$ cards and $2m-n$ players with $1$ card. Consider the game state after one turn. If a person with $2$ cards wins the hand, we will have $n-m-1$ players with $1$ card and one player with $m+1$ cards. If a player with $1$ card wins the hand, we will have $n-m$ players with $1$ card and one player with $m$ cards. If the player with most of the cards wins again with probability roughly $1/o(m) \simeq 0$ for large $o(m),$ then the game will have taken $2$ turns. Otherwise, if any of the other players wins, they will have around $n-m = o(m)$ cards or $n-m+1 = o(m)+1$ cards depending on which case while the remaining player will have $m$ or $m-1$ cards depending on the case. From this position, the game will thus take $m\cdot o(m)$ turns to complete, so the total number of turns is roughly $2+m\cdot o(m).$ By choosing larger functions $o(m)$ until reaching $cm$ for $c > 1$, we see that the total number of turns moves from $1$ when $n = m$ to $2+ m \cdot o(m)$ when $n = m + o(m)$ to $\Theta(n^2)$ when $n$ is $cm$ for $c > 1$ asymptotically.
\end{remark}

\section{Simulations}\label{sec:sim}

We ran simulations of the simple sticky random walk in order to analyze its termination time. We assume that $m$ divides $n$ and the $n$ cards are distributed equally among the players at the start of the turn. For each card count and player count pair, we simulated $500$ walks and determined the average termination time in order to estimate the constant of the $\Theta(n^2).$ The results are shown in \cref{table}.

\begin{figure}[h]
    \begin{subfigure}[b]{0.32\textwidth}
    \centering
    \caption{2 players}
    \includegraphics[width=\linewidth]{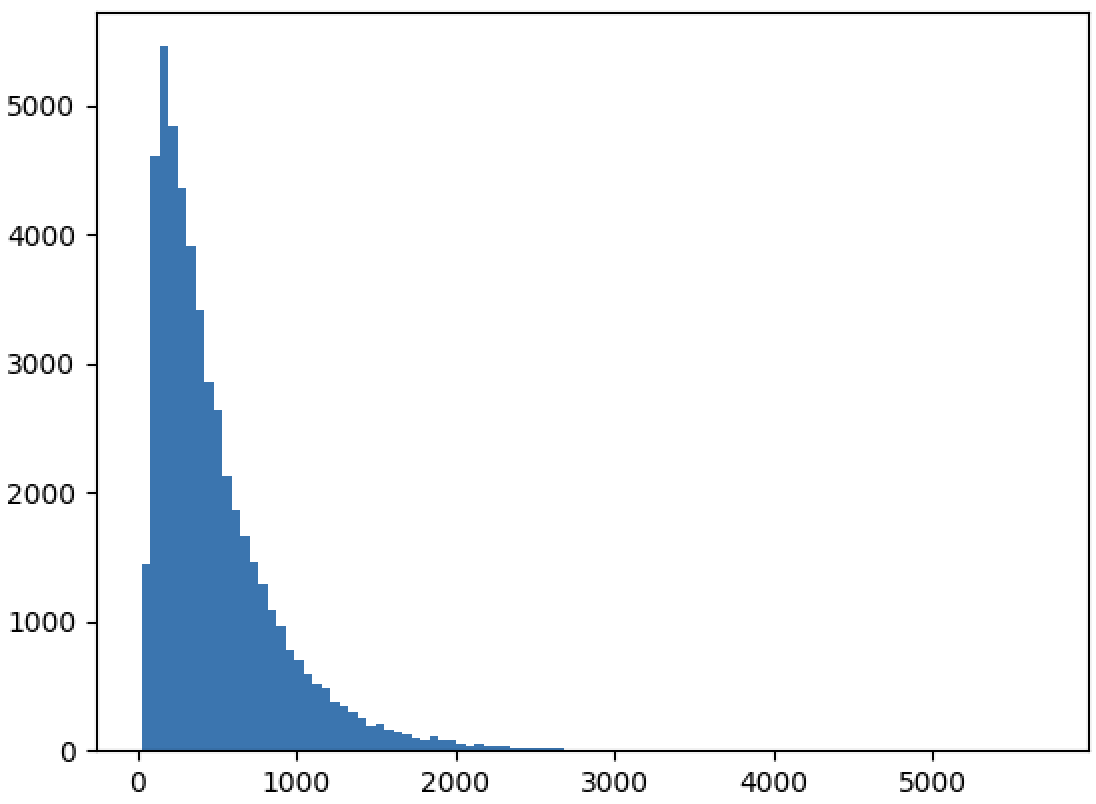}
    \\Mean \#rounds: $488$ \\Median \#rounds: $370$ \\Max \#rounds: $5692$
    \end{subfigure}
    \begin{subfigure}[b]{0.32\textwidth}
    \centering
    \caption{4 players}
    \includegraphics[width=\linewidth]{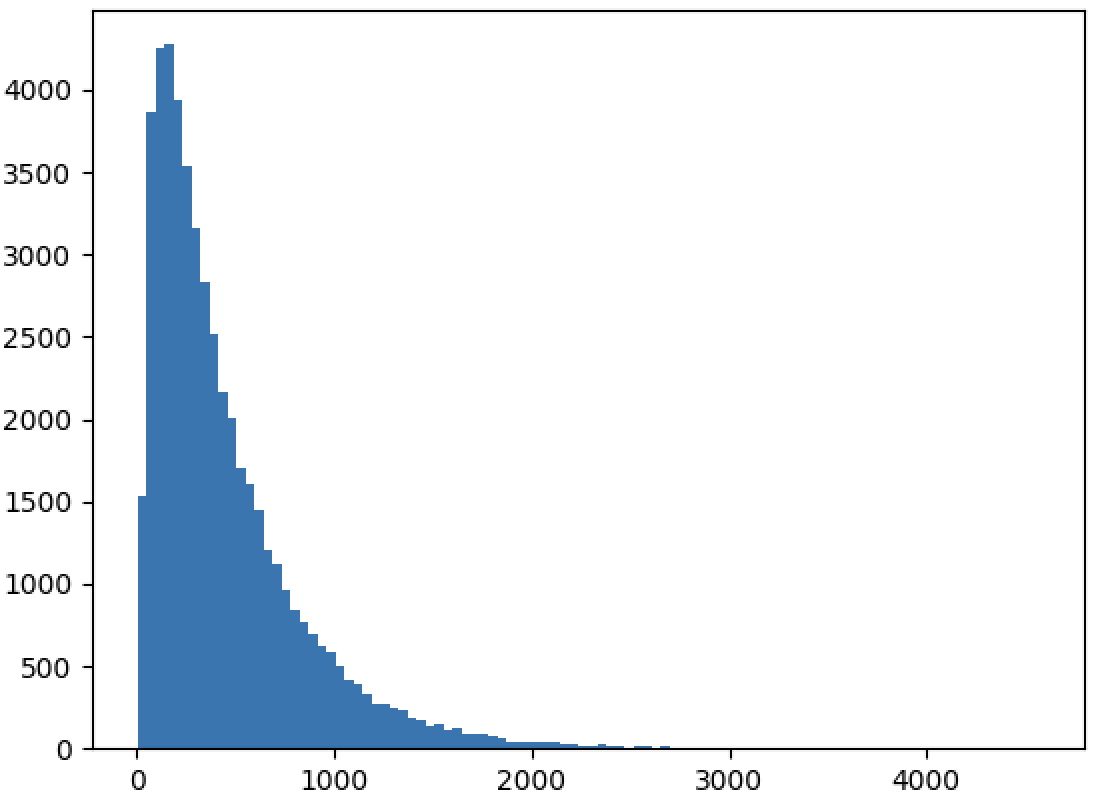}
    \\Mean \#rounds: $447$ \\Median \#rounds: $329$ \\Max \#rounds: $4569$
    \end{subfigure}
    \begin{subfigure}[b]{0.32\textwidth}
    \centering
    \caption{13 players}
    \includegraphics[width=\linewidth]{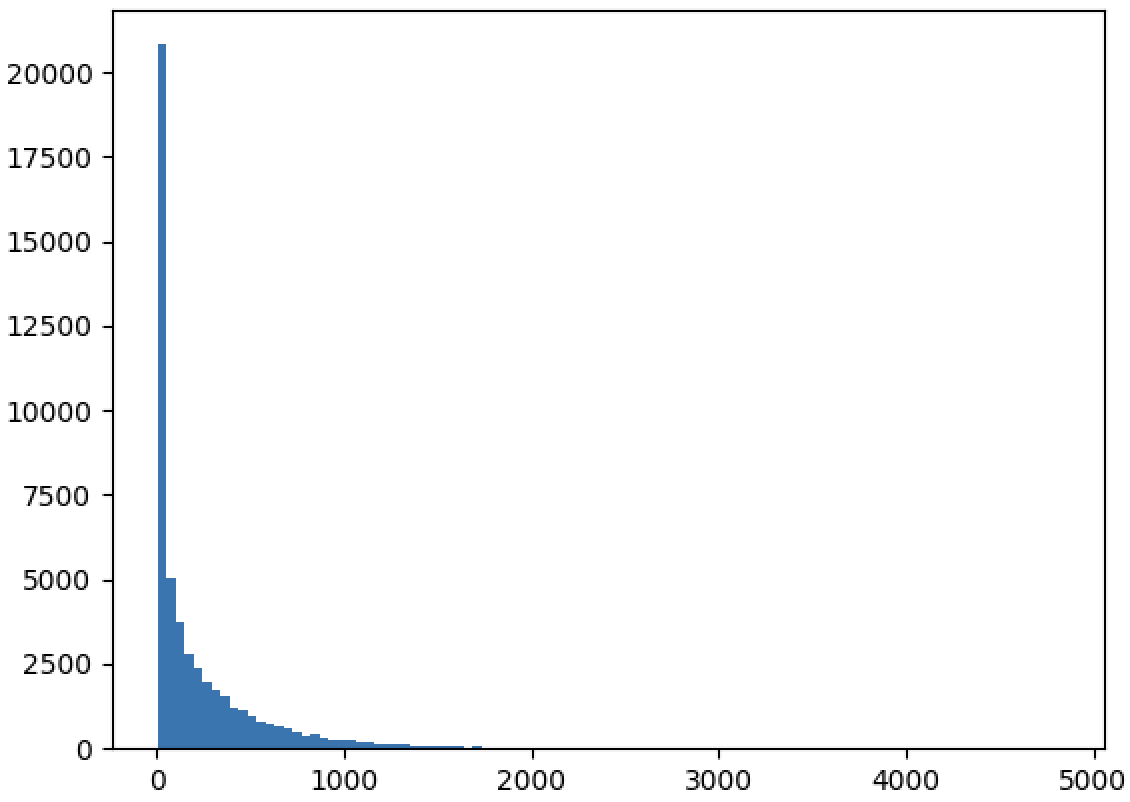}
    \\Mean \#rounds: $222$ \\Median \#rounds: $89$ \\Max \#rounds: $4815$
    \end{subfigure}
    \caption{Simulations of a real game of War with $2,$ $4,$ and $13$ players in the game.}
    \label{fig}
\end{figure}

\begin{table}[!htb]
\centering
\begin{tabular}{|c|c|c|c|}
\hline
$m$ & $n$ & $\mathbb{E}[\tau] \pm 2\cdot \text{SEM}$ & $\mathbb{E}[\tau]/n^2$ \\ 
\hline\hline
3 & 12 & $35.68 \pm 2.42$ & 0.248 \\ 
3 & 48 & $570.02 \pm 39.66$ & 0.247 \\ 
\hline
4 & 16 & $63.60 \pm 4.92$ & 0.248 \\ 
4 & 64 & $1128.33 \pm 78.18$ & 0.275 \\ 
\hline
5 & 20 & $92.78 \pm 7.22$ & 0.232 \\ 
5 & 80 & $1722.16 \pm 124.00$ & 0.269 \\ 
\hline
16 & 512 & $67311.82 \pm 4453.60$ & 0.257 \\ 
32 & 512 & $67988.85 \pm 4801.96$ & 0.259 \\
64 & 512 & $71081.94 \pm 4726.70$ & 0.271 \\ 
128 & 512 & $67381.37 \pm 4841.60$ & 0.257 \\
\hline
16 & 1024 & $273226.59 \pm 17407.10$ & 0.261 \\ 
32 & 1024 & $280040.75 \pm 20923.90$ & 0.267 \\ 
64 & 1024 & $278605.77 \pm 19548.66$ & 0.266 \\ 
128 & 1024 & $289191.59 \pm 19605.60$ & 0.276 \\ 
\hline
\end{tabular}
\caption{Simulation results for specified $m$ and $n$ pairs. All data points were formed by averaging $500$ absorption times of simulated random walks.}
\label{table}
\end{table}

Based on the results, it appears that the $\EE[\tau]/n^2$ is increasing slowly, approaching a constant independent of $n$ or $m$ which is around $0.3$.


We also ran simulations of the real game of War with a standard 52-card deck. Ties were broken by war rounds. Any player who runs out of cards during a war round is immediately out of the game. For each player count, we play $50,000$ games. The results of the simulations are shown in \cref{fig}. In these results, we can again see that the average number of rounds does not change significantly as the number of players increases except at $m = 13$ which could be caused by initial war rounds both removing players quickly and creating one player with most of the cards.

\section{Approximate Formula}\label{sec:approxformula}
In this section, we propose an approximate formula for the expected amount of time until at least one player loses all of their cards given the distribution of cards to the players at the start of the game. This is not the same as the expected time when the game ends. We start by motivating why finding such a formula is important.

\subsection{Motivation}\label{subsec:motivation} 
After \cref{sec:absorptiontime}, the main goal is to gain more precise estimates of the expected time when the game terminates. For the rest of this section, we only consider the simple sticky random walk as it models $\mfp$-war and top card $\mfp$-war. If $m = 3,$ we can actually determine an exact formula for the expected length of the game. We use a similar strategy to the proof of \cref{thm:srwabsorptiontime}, but we use a simpler martingale for this argument
$$X_t = \sum_{i = 1}^m A_{i,t}^2 - \sum_{s = 0}^{t-1} |C_s|(|C_s|-1).$$
Recall that $C_s$ is the set of players who still have cards at time $s.$ The proof that this is a martingale is analogous to the proof of \cref{lem:cubesmartingale} and is thus omitted. As before, we can now use the optional stopping theorem.

Let $\tau$ be the expected time when one player has won the game. Let $\tau_1$ denote the first time that some player has run out of cards. Using a formula originally proposed by Sandell \cite{RePEc:eee:stapro:v:7:y:1988:i:1:p:61-63}, we have
$$
\EE[\tau_1] = \frac{A_{1, 0}A_{2, 0}A_{3, 0}}{n - 2}.
$$
Thus, we have that
\begin{align*}
    \EE\left[\sum_{s = 0}^{\tau - 1}|C_s|(|C_s| - 1)\right] &= \EE\left[\sum_{s = 0}^{\tau_1 - 1}|C_s|(|C_s| - 1) + \sum_{s = \tau_1}^{\tau - 1}|C_s|(|C_s| - 1)\right], \\
    &= \EE\left[\sum_{s = 0}^{\tau_1 - 1}3(2) + \sum_{s = \tau_1}^{\tau - 1}2(1)\right], \\
    &= 6\EE[\tau_1] + 2\EE[\tau] - 2\EE[\tau_1], \\
    &= 2\EE[\tau] + 4\EE[\tau_1].
\end{align*}
Combining this with the optional stopping theorem, we have
$$
    n^2 - A_{1, 0}^2 - A_{2, 0}^2 - A_{3, 0}^2= \EE\left[\sum_{s = 0}^{\tau - 1}|C_s|(|C_s| - 1)\right] = 2\EE[\tau] + 4\EE[\tau_1],
$$
which implies

$$    \EE[\tau] = \frac{n^2 - A_{1, 0}^2 - A_{2, 0}^2 - A_{3, 0}^2 - 4\EE[\tau_1]}{2}
              = A_{1, 0}A_{2, 0} + A_{1, 0}A_{3, 0} + A_{2, 0}A_{3, 0} -\frac{2A_{1, 0}A_{2, 0}A_{3, 0}}{n - 2}.
$$
If we assume that $A_{1,0}=A_{2,0}=A_{3,0}=\frac{n}{3}$, then $\EE[\tau] = n^2 (7/27 + o(1)).$

Knowing $\EE[\tau_1]$ was critical to this proof and thus finding approximates for $\EE[\tau_1]$ when $m > 3$ is a useful endeavor. \cite{RePEc:eee:stapro:v:7:y:1988:i:1:p:61-63} notes that their method using a martingale likely does not work for $m > 3.$

\subsection{Formula}

Let $\x$ be the vector of cards that people have. Let $m\geq 3$ be the number of entries in this vector. Denote by $\v_i$ where $1\leq i \leq m$ the vector of length $m$ with $m-1$ in the $i$th entry and a $-1$ in all other entries. As always, we will use $y_i$ to denote the $i$th entry of a vector $\y.$ Additionally, we will use $(y_j)_i$ to denote the $i$th entry of a vector $\y_j.$

Let $k(\x)$ be
$$\prod_{j = 1}^m x_j - \frac{1}{m} \sum_{i = 1}^m \prod_{j = 1}^m (x_j + (v_i)_j).$$
Our approximate formula $f^*(\x)$ will be
$$f^*(\x) = \begin{cases}
    \frac{m(m-1)}{4m-6} \cdot \frac{\prod_{i = 1}^m x_i}{k(\x)} \text{ if all entries in $\x$ are greater than $0.$ } \\
    0 \text{ otherwise}
\end{cases}.$$
The reason for $\frac{m(m-1)}{4m-6}$ is to ensure that for $\x = \langle \ell + a_1, \cdots, \ell + a_m\rangle,$
$$\lim_{\ell\to \infty} f^*(\x) - \frac{1}{m} \sum_{i = 1}^m f^*(\x+\v_i) = 1.$$ To show this is a tedious but elementary calculation and it is thus omitted. Just like with any rational polynomial, the limit is the quotient between the leading coefficient of the numerator and the leading coefficient of the denominator assuming both numerator and denominator are the same degree, so we just cross multiply and simplify coefficients until we get the leading coefficient of both polynomials.

Note that $f^*(\x)$ when $m = 3$ agrees with the formula of \cite{RePEc:eee:stapro:v:7:y:1988:i:1:p:61-63}. In what follows, we will prove an important property of $f^*(\x)$ concerning the difference between $f^*(\x)$ and $\sum_{i = 1}^m f^*(\x+\v_i).$ For the actual expected value $E(\x),$ we have
$$E(\x) = 1 + \frac{1}{m} \sum_{ i= 1}^m E(\x+\v_i).$$ We will show
$f^*(\x) \geq \frac{1}{m} \sum_{i = 1}^m f^*(\x+\v_i).$
Ultimately, we would want to show that
$$f^*(\x) \approx 1 + \frac{1}{m} \sum_{i = 1}^m f^*(\x+\v_i).$$
The intuition for why $f^*(\x)$ is then a good approximation for $E(\x)$ is that whenever an entry of $\x$ is $0,$ $f^*(\x) = E(\x).$ This is analogous to an initial condition. If we then showed,
$f^*(\x) \approx 1 + \frac{1}{m} \sum_{ i = 1}^m f^*(\x+\v_i),$
then $f^*(\x)$ would have the same growth as $E(\x).$ Thus, $f^*(\x) \approx E(\x).$ We believe that examining the Hessian of $f^*(\x)$ which is quite easy because of the structure of $k(\x)$ is useful in achieving this goal. 

The remainder of this section will be devoted to showing
$$f^*(\x) \geq \frac{1}{m} \sum_{i = 1}^m f^*(\x+\v_i).$$
It suffices to show that 
$$f(\x) \geq \frac{1}{m} \sum_{i = 1}^m f(\x+\v_i),$$
where $f(\x) = f^*(x) \cdot \frac{4m-6}{m(m-1)}.$

Without loss of generality, assume that the entries in $\x$ are ordered in non-decreasing order. We start by showing some general facts about $k(\x)$ and $f(\x).$ We have that $k(\x)$ is defined by
\begin{align*}
k(\x) &=  \prod_{i = 1}^m x_i - \frac{1}{m} \sum_{i = 1}^m \prod_{j = 1}^m (x_i + (v_i)_j),\\
&= \prod_{j = 1}^m x_j - \frac{1}{m} \sum_{i = 1}^m \left(\prod_{j = 1}^m (x_j-1)+ m\prod_{\substack{j = 1\\j\neq i}}^m (x_j-1)\right),\\
&= \prod_{j = 1}^m x_j - \prod_{j = 1}^m (x_j-1) - \sum_{i = 1}^m \prod_{\substack{j = 1\\j\neq i}}^m (x_j-1).
\end{align*}
We refer to the bottommost formula for $k(\x)$ as the \textit{product form} of $k(\x).$ As a result of the product form, if $\x$ has two entries which are $1,$ then the second and third terms must be $0,$ so $k(\x)$ is just the product of the $x_j,$ so $f(\x) = 1.$ Let $S_i$ denote the $i$th symmetric polynomial, i.e., the sum of the $i$th-wise products of the $x_j.$ Therefore,
\begin{align*}
k(\x) &= \prod_{j = 1}^m x_j - \prod_{j = 1}^m (x_j-1) - \sum_{i = 1}^m \prod_{\substack{j = 1\\j\neq i}}^m (x_j-1),\\
&= S_m - \sum_{i = 0}^m (-1)^{m-i}S_i - \sum_{i = 0}^{m-1} (-1)^{m-1-i} (m-i) S_i,\\
&= \sum_{i = 0}^{m-2} (-1)^{m-i} (m-1-i) S_{i}.
\end{align*}
We refer to the bottommost formula for $k(\x)$ as the \textit{sum form} of $k(\x).$ This product and sum form of $k(\x)$ will be important throughout the proofs. Because $k(\x)$ is linear with respect to each $x_i$ we then have $\frac{\partial^2}{\partial x_i^2} k(\x) = 0$ for all $i.$ We now illustrate one important property of $k(\x).$ Note that we assume $m\geq 2$ for all of these results.
\begin{lemma}
    For all $\x$ with all entries greater than or equal to $1,$ we have $k(\x) > 0$ and $\frac{\partial}{\partial x_i} k(\x) \geq 0$ for all $1\leq i \leq m.$ In addition $\frac{\partial}{\partial x_a} k(\x) \geq \frac{\partial}{\partial x_b} k(\x)$ if $a < b.$
\end{lemma}
\begin{proof}
    To show this, note that $k(\1) = 1$ where $\1$ is the all ones vector using the product formula for $k(\x).$ We will now show that $k(\x) > 0$ for all $i$ through induction on $m.$ When $m = 2,$ we have that $k(\x) = 1 > 0.$ Now assume that $k(\x) > 0$ when the number of entries in $\x$ is $m-1.$ Let the number of entries in $\x$ be $m$ now. Using the product form of $k(\x),$ we have
    $$\frac{\partial}{\partial x_i} k(\x) = \prod_{\substack{j = 1\\j \neq i}}^m x_j - \prod_{\substack{j = 1\\j \neq i}}^m (x_j-1) - \sum_{\substack{\ell = 1 \\ \ell \neq i}}^m \prod_{\substack{j = 1\\j \not \in \{i,\ell\}}}^m (x_j - 1) = k(\x') \geq 0,$$
    where $\x'$ denotes the vector $\x$ without the $x_i$ entry. Thus, $k(\x)$ increases if any of its entries increase. Because $k(\1) > 0,$ we must have $k(\x) > 0$ for all $\x$ then. In addition because $\frac{\partial}{\partial x_i} k(\x) = k(\x') > 0$ for $m\geq 3$ and when $m = 2$ we have that $k(\x)$ is constant so  $\frac{\partial}{\partial x_i} k(\x)$ is $0,$ we must have that $\frac{\partial}{\partial x_i} k(\x) = k(\x') \geq 0$ for $m\geq 2.$

    We have showed that $k(\x)$ increases when any of the components of $\x$ increase, and we also showed that $\frac{\partial}{\partial x_i} k(\x) = k(\x')$ where $\x'$ is $\x$ without the $x_i$ entry, so $k(\x')$ is the largest when we remove the smallest entry, so $\frac{\partial}{\partial x_a} k(\x) \geq \frac{\partial}{\partial x_b} k(\x)$ if $a < b.$
\end{proof}

We are ready to consider the case when $x_1 = 1.$ In this case, because $f(\x+\v_i) = 0$ for $i\geq 2$, we just need to show.
\begin{lemma}
If $x_1 = 1,$
$$f(\x) \geq \frac{1}{m} f(\x+\v_1).$$
\end{lemma}
\begin{proof}

    Let $\h_i$ be the vector with a $1$ in the first entry, a $-1$ in the $i$th entry, and zeros for all other entries. Order the entries of $\x$ in non-decreasing order. Therefore, $x_1 = 1.$ Let $2\leq j \leq m$ be an integer. Consider the relationship between $f(\y)/(j-1)$ and $f(\y+\h_j)/j$ where $\y$ is the sum of $\x$ and $\h_i$ for $2\leq i \leq j-1.$ Notice that
    $$\frac{f(\y)}{j-1} = \frac{1}{j-1} \cdot \frac{1}{k(\y)} \prod_{i = 1}^m y_i = \frac{1}{k(\y)} \prod_{i = 2}^m y_i,$$
    and
    $$\frac{f(\y+\h_j)}{j} = \frac{1}{j} \cdot \frac{1}{k(\y+\h_j)} \prod_{i = 1}^m (y_i+(h_j)_i) = \frac{1}{k(\y+\h_j)} \frac{y_j-1}{y_j} \prod_{i = 2}^m y_i.$$
    Thus, if we show that $y_j k(\y+\h_j) \geq (y_j-1) k(\y),$ then this implies that $f(\y)/(j-1) \geq f(\y+\h_j)/j,$ so $f(\y) \geq \frac{j-1}{j} f(\y+\h_j).$ This would imply
    $$f(\x) \geq \frac{1}{2} f(\x+\h_2) \geq \frac{1}{3} f(\x+\h_2+\h_3) \geq \cdots \geq \frac{1}{m} f(\x+\h_2+\cdots+\h_m) = \frac{1}{m} f(\x+\v_i).$$

    To show $y_j k(\y+\h_j) \geq (y_j-1) k(\y),$ we use induction on $m.$
    If $m = 2,$ then $k(\x) = 1,$ so $k(\y+\h_j) = k(\y),$ so $y_j k(\y+\h_j) \geq (y_j-1) k(\y).$ If $m = 3,$ then $k(\x) = x_1 + x_2 + x_3 - 2$ which is also constant, so $k(\y+\h_j) = k(\y),$ so $y_j k(\y+\h_j) \geq (y_j-1) k(\y)$ as well.
    
    Assume that the inequality holds for when $\y$ is of size $m-1.$ Consider the case when $y_i = 1$ for $2\leq i \leq m$ where $i \neq j.$ For $m\geq 4$ then for both $\y$ and $\y+\h_j,$ they both have at least two entries which are $1.$ Therefore, using the product form of $k(\y+\h_j),$
    $$y_j k(\y+\h_j) = y_j \prod_{i = 1}^m (y_i + (h_j)_i)  = x_j\cdot j(y_j-1) \geq (y_j-1) (j-1) y_j = (y_j-1) \prod_{i = 1}^m y_i = (y_j-1) k(\y).$$
    Now let us show that increasing $y_i$ for $2\leq i \leq m$ where $i\neq j$ makes the difference between $y_j k(\y+\h_j)$ and $(y_j-1) k(\y)$ larger. It suffices to show that
    $$y_j \frac{\partial}{\partial y_i} k(\y+\h_j) \geq (y_j-1) \frac{\partial}{\partial y_i} k(\y),$$
    for $2\leq i \leq m$ where $i\neq j.$ Recall that $\frac{\partial}{\partial x_i} k(\x)$ is $k(\x')$ where $\x'$ is $\x$ without the $x_i$ entry, so this is true by the inductive hypothesis.
\end{proof}

We now consider the case where the first entry is larger than $1.$ Therefore all entries are larger than $1$ as they are in nondecreasing order. The way we will prove the inequality is we will show the following set of inequalities
\begin{align*}
    f(\x) &\geq \frac{1}{m} f(\x+\v_1) + \frac{m-1}{m} f\left(\x-\frac{1}{m-1} \v_1\right).\\
    f\left(\x-\frac{1}{m-1} \v_1\right) &\geq \frac{1}{m-1} f(\x+\v_2) + \frac{m-2}{m-1} f\left(\x-\frac{1}{m-2} \v_2 -\frac{1}{m-2} \v_1 \right).\\
    &\,\,\, \vdots \\
    f\left( \x - \frac{1}{2}\sum_{i = 1}^{m-2} \v_i \right) &\geq \frac{1}{2} f(\x+\v_{m-1}) + \frac{1}{2} f\left(\x-\sum_{i = 1}^{m-1} \v_{m-1}\right) = \frac{1}{2}f(\x+\v_{m-1}) + \frac{1}{2} f(\x+\v_m).
\end{align*}

Notice that we can rewrite these inequalities using the vectors $\u_i$ where $\u_i$ has an $m-i$ in the $i$th entry, zeros in all entries before the $i$th entry, and a $-1$ in all the entries after the $i$th entry.
\begin{align*}
    f(\x) &\geq \frac{1}{m} f(\x+\v_1) + \frac{m-1}{m} f\left(\x-\frac{1}{m-1} \u_1\right).\\
    f\left(\x-\frac{\u_1}{m-1}\right) &\geq \frac{1}{m-1} f\left(\x+\v_2\right) + \frac{m-2}{m-1} f\left(\x-\frac{\u_1}{m-1} - \frac{m \u_2}{(m-1)(m-2)} \right).\\
    &\, \, \, \vdots \\
    f\left( \x - \sum_{i = 1}^{m-2} \frac{m\u_i}{(m-i+1)(m-i)} \right) &\geq \frac{1}{2} f(\x+\v_{m-1}) + \frac{1}{2} f\left(\x-\sum_{i = 1}^{m-1} \frac{m\u_i}{(m-i+1)(m-i)}\right)
\end{align*}

We will prove each of these inequalities through Jensen's inequality about concavity. Let us first show another important property of $k(\x)$ first though.

\begin{lemma}
    For any $\x$ with all entries greater than or equal to $1,$ we have
    $$k(\x) - x_i \frac{\partial}{\partial x_i} k(\x) \geq 0.$$
    Furthermore,
    $$k(\x) - x_a \frac{\partial}{\partial x_a} k(\x) \geq k(\x) - x_b \frac{\partial}{\partial x_b} k(\x),$$
    if $a<b.$
\end{lemma}
\begin{proof}
    Notice that $k(\x) - x_i \frac{\partial}{\partial x_i} k(\x)$ are all of the terms in $k(\x)$ without $x_i.$ Using the sum form of $k(\x)$ we then have
    $$k(\x) - x_i \frac{\partial}{\partial x_i} k(\x) = \sum_{j = 0}^{m-2} (-1)^{m-j}(m-1-j) T_j,$$
    where $T_j$ is the $j$th-wise product of $x_1$ through $x_m$ excluding $x_i.$ Consider the derivative of
    $$\prod_{j = 1, j \neq i} (x_j-x) = \sum_{j = 0}^{m-1} (-1)^{m-1-j} T_{j} x^{m-1-j},$$
    which is
    $$-\sum_{\substack{\ell = 1\\ \ell \neq i}}^m \prod_{\substack{j = 1\\ j\not \in \{i,\ell\}}} (x_j-x) = \sum_{j = 0}^{m-2} (-1)^{m-1-j} (m-1-j)T_j x^{m-2-j}.$$
    Negating both sides and evaluating this derivative at $x = 1,$ we then have 
    $$\sum_{\substack{\ell = 1\\ \ell \neq i}}^m \prod_{\substack{j = 1\\ j\not \in \{i,\ell\}}} (x_j-1) = \sum_{j = 1}^{m-1} (-1)^{m-j} (m-1-j)T_j = k(\x) - x_i \frac{\partial}{\partial x_i} k(\x).$$
    Notice that the products in the leftmost expression are nonnegative as $x_j > 1$ for all $j,$ so $k(\x) - x_i \frac{\partial}{\partial x_i} k(\x) > 0.$ Furthermore, if $1\leq a < b \leq m$ are integers, then
    \begin{align*}
        k(\x) - x_a \frac{\partial}{\partial x_a} k(\x) &= \sum_{\substack{\ell = 1\\ \ell \neq a}}^m \prod_{\substack{j = 1\\ j\not \in \{a,\ell\}}} (x_j-1)\\
        &= \sum_{\substack{\ell = 1 \\ \ell \not \in \{a,b\}}}^m \prod_{\substack{j = 1 \\ j \not \in \{a,\ell\}}}^m (x_j-1) + \prod_{\substack{j = 1 \\ j \not \in \{a,b\}}}^m (x_j-1), \\
        &\geq \sum_{\substack{\ell = 1 \\ \ell \not \in \{a,b\}}}^m \prod_{\substack{j = 1 \\ j \not \in \{b,\ell\}}}^m (x_j-1) + \prod_{\substack{j = 1 \\ j \not \in \{a,b\}}}^m (x_j-1),\\
        &= \sum_{\substack{\ell = 1\\ \ell \neq b}}^m \prod_{\substack{j = 1\\ j\not \in \{b,\ell\}}} (x_j-1) = k(\x) - x_b \frac{\partial}{\partial x_b} k(\x). \\
    \end{align*}
\end{proof}

Using this lemma, we can now bound the second partial derivatives of $f.$

\begin{lemma}
    We have
    $$\frac{\partial^2}{\partial x_i^2} f(\x) \leq 0,$$
    for all $1\leq i \leq m$ if all the entries in $\x$ are greater than or equal to $1.$
\end{lemma}
\begin{proof}
    Notice that
    $$\frac{\partial}{\partial x_i} f(\x) = \frac{\partial}{\partial x_i} \frac{\prod_{j = 1}^m x_j}{k(\x)} = \frac{\prod_{\substack{j = 1, j\neq i}}^m x_j}{k(\x)} - \frac{\left(\frac{\partial}{\partial x_i} k(\x) \right)\prod_{j = 1}^m x_j}{k(\x)^2},$$
    so
    \begin{align*}
    \frac{\partial^2}{\partial x_i^2} f(\x) &= -2 \frac{\left(\frac{\partial}{\partial x_i} k(\x) \right)\prod_{j = 1, j \neq i}^m x_j}{k(\x)^2} + 2\frac{\left(\frac{\partial}{\partial x_i} k(\x) \right)^2\prod_{j = 1}^m x_j}{k(\x)^3},\\
    &= 2\frac{\left(\frac{\partial}{\partial x_i} k(\x) \right)\prod_{j = 1, j \neq i}^m x_j}{k(\x)^3} \left(-k(\x) + x_i \frac{\partial}{\partial x_i} k(\x) \right) \leq 0.
    \end{align*}
\end{proof}

\begin{lemma}
    We have
    $$\frac{\partial^2}{\partial x_b \partial x_a} f(\x) \geq 0,$$
    and
    $$\frac{\partial^2}{\partial x_b \partial x_a} f(\x) \geq \frac{\partial^2}{\partial x_c \partial x_a} f(\x),$$
    if $b < c$ for all $\x$ with entries greater than or equal to $1.$
\end{lemma}
\begin{proof}
    We start by evaluating the second derivative of $f(\x)$ with respect to $x_a$ and $x_b.$
    We have
    $$\frac{\partial}{\partial x_a} f(\x) = \frac{\prod_{\substack{j = 1, j\neq a}}^m x_j}{k(\x)} - \frac{\left(\frac{\partial}{\partial x_a} k(\x) \right)\prod_{j = 1}^m x_j}{k(\x)^2},$$
    and we have that $\frac{\partial^2}{\partial x_b \partial x_a} f(\x)$ is 
    \begin{align*}
    &\frac{\prod_{\substack{j = 1, j\not \in \{a,b\}}}^m x_j}{k(\x)} - \frac{\left(\frac{\partial}{\partial x_b} k(\x) \right)\prod_{\substack{j = 1 \\ j \neq a}}^m x_j}{k(\x)^2} - \frac{\left(\frac{\partial}{\partial x_a} k(\x) \right)\prod_{\substack{j = 1 \\ j \neq b}}^m x_j}{k(\x)^2} + 2 \frac{\left(\frac{\partial}{\partial x_b} k(\x) \right)\left(\frac{\partial}{\partial x_a} k(\x) \right)\prod_{j = 1}^m x_j}{k(\x)^3},\\
    = \, &\frac{\prod_{\substack{j = 1, j\not \in \{a,b\}}}^m x_j}{k(\x)^3}\left( k(\x) - x_a \frac{\partial}{\partial x_a} k(\x)\right) \left( k(\x) - x_b \frac{\partial}{\partial x_b} k(\x)\right) +  \frac{\left(\frac{\partial}{\partial x_b} k(\x) \right)\left(\frac{\partial}{\partial x_a} k(\x) \right)\prod_{j = 1}^m x_j}{k(\x)^3} > 0.
    \end{align*}
    Because $\prod_{\substack{j = 1, j\not \in \{a,b\}}}^m x_j$ is larger when $b$ is smaller as we are removing the smallest element and $k(\x) - \frac{\partial}{\partial x_i} k(\x),$ is larger for smaller $i$ and $\frac{\partial}{\partial x_i} k(\x)$ is larger for smaller $i,$ meaning that nonconstant portion of every part of the expression for $\frac{\partial^2}{\partial x_b \partial x_a} f(\x)$ increases for smaller $i,$ we also see that 
    $$\frac{\partial^2}{\partial x_b \partial x_a} f(\x) \geq \frac{\partial^2}{\partial x_c \partial x_a} f(\x).$$
\end{proof}

\begin{theorem}
    For all $\x$ with entries greater than or equal to $1,$ we have
    $$f(\x) \geq \frac{1}{m} \sum_{i = 1}^m f(\x+\v_i).$$
\end{theorem}
\begin{proof}

    We have already shown this inequality in the case when $\x$ has an entry equal to $1.$ Therefore, it suffices to prove the inequality in the case that all entries are greater than or equal to $2.$ Recall the inequalities using the $\u_i$ which are equivalent to the theorem statement. To start, let us show that all of the vectors in the inequalities have components with entries greater than or equal to $1.$ First $\x+\v_i$ subtracts at most $1$ from each of the components, so each component must be greater than or equal to $1$ as each component of $\x$ is greater than or equal to $2.$ Furthermore for any of the first $j$ entries of
    $$\x-\sum_{i = 1}^j \frac{m\u_i}{(m-i+1)(m-i)},$$
    we see they are preserved after by subtracting
    $$\x-\sum_{i = 1}^j \frac{m\u_i}{(m-i+1)(m-i)}- \sum_{i = j+1}^{m-1} \frac{m\u_i}{(m-i+1)(m-i)} = \x + \v_m,$$
    so each of the first $j$ entries is the same as the first $j$ entries of $\x+\v_m$ which are all greater than or equal to $1.$ Any of the entries after the first $j$ only increase as the entries after the first $j$ for each of the $\u_i$ for $1\leq i \leq j$ are all negative, so subtracting them from $\x$ only increases the entries. Therefore, they are all greater than or equal to $1.$

    Let $H$ denote the Hessian of $f,$ i.e.,
    $$H(\x) = \begin{bmatrix}
        \frac{\partial^2}{\partial x_1^2} f(\x) & \cdots & \frac{\partial^2}{\partial x_1 \partial x_n} f(\x) \\
        \vdots & \ddots & \vdots \\
        \frac{\partial^2}{\partial x_n \partial x_1} f(\x) & \cdots & \frac{\partial^2}{\partial x_n^2} f(\x)
    \end{bmatrix}.$$
    Let us show that $\u_i^T H(\x) \u_i \leq 0$ for all $1\leq i \leq m-1$ and and $\x$ with entries greater than or equal to $1.$ This would imply all of the inequalities because we would then know that $f$ is concave along each of the $\u_i$ for $\x$ with entries greater than or equal to $1$ which implies the inequalities by Jensen's inequality.

    We have
    \begin{align*}
        \u_i^T H(\x) \u_i &= (m-i)^2 \frac{\partial^2}{\partial x_i^2} f(\x) - 2(m-i) \sum_{j = i+1}^m \frac{\partial^2}{\partial x_j \partial x_i} f(\x) + \sum_{j = i+1}^m \sum_{\ell = i+1}^m \frac{\partial^2}{\partial x_j \partial x_\ell} f(\x),\\
        &\leq -(m-i-1) \sum_{j = i+1}^m \frac{\partial^2}{\partial x_j \partial x_i} f(\x) + \sum_{j = i+1}^m \sum_{\substack{\ell = i+1 \\ \ell \neq j}}^m \frac{\partial^2}{\partial x_j \partial x_\ell} f(\x),\\
        &= \sum_{j = i+1}^m \sum_{\substack{\ell = i+1 \\ \ell \neq j}}^m \left(\frac{\partial^2}{\partial x_j \partial x_\ell} f(\x) - \frac{\partial^2}{\partial x_j \partial x_i} f(\x)\right),\\
        &\leq 0.
    \end{align*}
\end{proof}

\section{Further Directions}\label{sec:furtherdirec} 
As mentioned in \cref{subsec:motivation}, the main direction of interest is to find more precise estimates for the expected termination time of a game of War, specifically $\mfp$-war which boils down to a simple sticky random walk. In \cref{subsec:motivation}, we show a potential way of doing this which requires knowing the expected amount of time until some player loses their cards which is an interesting question in and of itself. A good step forward would be to solve the following problem.
\begin{problem}\label{prob}
    Find a relationship (e.g. an inequality) between $f^*(\x),$ the approximate formula in \cref{sec:approxformula}, and $E(\x),$ the true expected time until some player loses their cards.
\end{problem}
Furthermore, when $m\geq 4,$ there is an additional complication in the argument of \cref{subsec:motivation} in that we must know the expected distribution of the remaining cards after some number of people have lost all of their cards. For example, in the $m = 4$ case, we would need the expected time until at least one person has lost all of their cards and the expected time between at least one person losing all of their cards and at least two people losing their cards which involves knowing the distribution of cards after at least one person has lost all of their cards. An upper bound on $E(\x)$ is likely to be concave, so using Jensen's inequality eliminates the need to know the complete distribution. This would yield a lower bound on $\EE[\tau]$ where $\tau$ is the expected time when the game terminates. More work is needed to lower-bound $E(\x).$

\section*{Acknowledgements}
\noindent AA and NK were supported by MIT UROP+. EM was supported in part by a Simons Investigator Award, Vannevar Bush Faculty Fellowship ONR-N00014-20-1-2826, and ARO MURI W911NF1910217.

\bibliographystyle{alpha}
\bibliography{citations}






\end{document}